\documentclass[twoside]{amsart}
\usepackage{amsmath,amsthm,amssymb}

\theoremstyle{plain}
\newtheorem{theorem}{Theorem}[section]
\newtheorem{proposition}[theorem]{Proposition}
\newtheorem{lemma}[theorem]{Lemma}

\newtheorem{corollary}[theorem]{Corollary}

\theoremstyle{definition}
\newtheorem{remark}[theorem]{Remark}
\theoremstyle{remark}
\newtheorem{example}[theorem]{Example}
\newtheorem{acknowledgments}{Acknowledgments}

\newcommand{\A}{{\mathbb A}}
\newcommand{\C}{{\mathbb C}}
\newcommand{\F}{{\mathbb F}}

\newcommand{\I}{{\mathbb I}}
\newcommand{\M}{{\mathbb M}}

\newcommand{\Q}{{\mathbb Q}}
\newcommand{\R}{{\mathbb R}}

\newcommand{\Z}{{\mathbb Z}}
\newcommand{\llceil}{\left\lceil}
\newcommand{\rrceil}{\right\rceil}

\def\bar{\overline}

\newcommand{\cL}{ {\mathcal L} }
\newcommand{\cM}{ {\mathcal M} }
\newcommand{\cO}{ {\mathcal O} }
\newcommand{\cP}{ {\mathcal P} }
\newcommand{\cU}{ {\mathcal U} }
\newcommand{\cW}{ {\mathcal W} }
\newcommand{\cX}{ {\mathcal X} }
\newcommand{\cY}{ {\mathcal Y} }

\newcommand{\ord}{ {\rm ord} }
\newcommand{\sgn}{ {\rm sgn} }
\newcommand{\NP }{ {\rm  NP} }     %Newton polygon
\newcommand{\HP }{ {\rm  HP} }     %Hodge polygon
\newcommand{\GNP }{ {\rm  GNP} }     %Newton polygon
\newcommand{\Tr}{ {\rm Tr} }

\begin{document}

%Topmatter
\title{Asymptotic
variation of $L$ functions of one-variable exponential sums }
\author{Hui June Zhu}
\address{
Department of mathematics and statistics,
McMaster University,
Hamilton, ON L8S 4K1,
CANADA.
}
\email{zhu@cal.berkeley.edu}
\date{Sept. 13, 2003}
\subjclass{11,14}

\begin{abstract}
Fix an integer $d\geq 3$.
Let $\A^d$ be the
dimension-$d$ affine space over the algebraic closure $\bar\Q$
of $\Q$, identified with the
coefficient space of degree-$d$ monic polynomials $f(x)$ in
one variable $x$.
For any $f(x)$ in $\A^d(\bar\Q)$,
let $\Q(f)$ be the field generated by coefficients of $f$
in $\bar\Q$. For each prime $p$ coprime to $d$,
pick an embedding from $\bar\Q$ to $\bar\Q_p$,
once and for all.
Let $\cP$ be the place in $\Q(f)$ lying over $p$
specified by the embedding of $\bar\Q$ in $\bar\Q_p$
(with residue field $\F_q$ say).
Suppose $f \in \A^d(\bar\Q\cap \bar\Z_p)$,
let $\NP(f(x)\bmod\cP)$ denote the $q$-adic Newton polygon of
the $L$ function $L(f(x)\bmod\cP;T)$ of exponential sums
of $f\bmod\cP$.
We prove that
there is a Zariski dense open subset $\cU$ defined over $\Q$ in
$\A^d$ such that for every geometric point $f(x)$ in $\cU(\bar\Q)$
and $p$ large enough (depending only on $f$)
one has
$\NP(f\bmod \cP)=\GNP(\A^d;\F_p)$ and
$$\lim_{p\rightarrow\infty}\NP(f(x)\bmod\cP)=\HP(\A^d),$$
where $\GNP(\A^d;\F_p)$ and $\HP(\A^d)$
are the generic Newton polygon and
the Hodge polygon, respectively (see \cite{Zhu:1}).
\end{abstract}

\maketitle

\section{Introduction}\label{S:Introduction}

In this paper we fix an integer $d\geq 3$. Let $\bar\Q$ be
the algebraic closure of $\Q$.
Let $\A^d$ be
the affine variety of dimension $d$ over $\bar\Q$, identified with
the coefficient space of degree-$d$ monic
polynomials $f(x)$ in one variable
$x$.
For any $f\in\A^d(\bar\Q)$ let
$\Q(f)$ denote the field generated by coefficients
of $f$ over $\Q$.
Let $p$ be any prime coprime to $d$.
Let $\bar\Q_p$ be the algebraic closure of $\Q_p$
and let $\bar\Z_p$ be its ring of integers.
For each $p$ we pick
an embedding from $\bar\Q$
into $\bar\Q_p$.
Let $\cP$ be the place in  $\Q(f)$ lying over $p$
specified by the embedding.
Henceforth we tacitly understand that
such embeddings are already picked once and for all.
Suppose the residue field at $\cP$ is $\F_q$ for
$q=p^a$ for some $a$.
Let $\ord_p(\cdot)$ denote
the $p$-adic valuation in an extension of $\Q_p$ with
$\ord_p(p)=1$; let $\ord_q(\cdot)$ be the $q$-adic valuation, i.e.,
$\ord_q(\cdot):=\frac{1}{a}\ord_p(\cdot)$. Let
$E(x)=\exp(\sum_{j=0}^{\infty}\frac{x^{p^j}}{p^j})$ be the
Artin-Hasse $p$-adic exponential function.
Let $\gamma$ be a root of $\log E(x)$ in
$\bar\Q_p$ with $\ord_p (\gamma)=\frac{1}{p-1}$. Then $E(\gamma)$
is a primitive $p$-th root of unity. We fix this $p$-th root of
unity for the entire paper and denote it by $\zeta_p$.  Note that
$\Z_p[\gamma]=\Z_p[\zeta_p]$.

Let $f\in(\bar\Z_p\cap\bar\Q)[x]$ be a degree-$d$ monic polynomial.
For every positive integer $\ell$ let
\begin{eqnarray}\label{E:S_ell}
S_\ell(f\bmod\cP)
:=\sum_{x\in\F_{q^\ell}}\zeta_p^{\Tr_{\F_{q^\ell}/\F_p}(f(x)\bmod\cP)}.
\end{eqnarray}
Then the $L$ function of the exponential sum of $f$ over $\F_q$ is defined by
\begin{eqnarray}\label{E:Lfunction}
L(f\bmod\cP;T)
:=\exp(\sum_{\ell=1}^{\infty}S_\ell(f\bmod\cP)\frac{T^\ell}{\ell}).
\end{eqnarray}
It is well known
that $L(f\bmod\cP;T)$ is a polynomial
in $1+T\Z[\zeta_p][T]$ of degree $d-1$
(e.g., see remarks in the Introduction of \cite{Zhu:1}).
So we may write,
\begin{eqnarray}\label{E:20}
L(f\bmod\cP;T)&=&1+b_1(f)T+b_2(f)T^2+\ldots+b_{d-1}(f)T^{d-1}
\end{eqnarray}
for some $b_n(f)\in\Z[\zeta_p][T]$. It is easy to see that
$L(f\bmod\cP;T)$ becomes independent of the choice of embedding
$\bar\Q\hookrightarrow\bar\Q_p$ if $p$ is large enough.

For any polynomial $\sum_{n=0}^{k}c_nT^n$ over $\bar\Q$ let
$\NP_q(\sum_{n=0}^{k}c_nT^n)$ denote its $q$-adic Newton polygon,
i.e., the lower convex hull in $\R^2$ of the points
$(n,\ord_q(c_n))$ with $0\leq n\leq k$. Now let
$\NP(f\bmod\cP):=\NP_q(L(f\bmod\cP;T))$. We
shall note below that $\NP(f\bmod\cP)$ is independent of the
choice of embedding $\bar\Q\hookrightarrow\bar\Q_p$.
By the Dieudonn\'e-Manin classification (see
\cite{Manin:1}), the Newton polygon of an abelian variety over a
finite field is determined by certain `formal isogeny types'.
Hence the Newton polygon of a smooth projective
curve over a finite field $\F_q$, same as that of its Jacobian
variety, is independent of the choice of $\F_q$.
Let $\NP(X_f\bmod\cP)$ be the
Newton polygon of the Artin-Schreier curve
$X_f$ given by affine equation $y^p-y=f\bmod\cP$.
One knows that $\NP(f\bmod\cP)=\NP(X_f\bmod\cP)/(p-1)$ where the latter Newton
polygon is shrunk by a factor of $p-1$ horizontally and
vertically (see \cite[Introduction]{Zhu:1}), all these above imply that
$\NP(f\bmod\cP)$ is independent of the
choice of embedding of $\bar\Q\hookrightarrow\bar\Q_p$.

The {\em Hodge polygon} $\HP(\A^d)$ of $\A^d$ is the
lower convex hull in $\R^2$ of the points $(n,\frac{n(n+1)}{2d})$
with $0\leq n\leq
d-1$. It is known that $\HP(\A^d)$ is a lower bound of $\NP(f\bmod\cP)$
(see \cite[Propositions 2.2 and 2.3]{Wan:1})
and that for every $f\in\A^d(\bar\Q\cap\bar\Z_p)$
one has $\NP(f\bmod\cP)=\HP(\A^d)$
if and only if $p\equiv 1\bmod d$
(see \cite[(3.11)]{Adolphson-Sperber}).
Our Hodge polygon, which inherits that from \cite{Wan:1, Wan:2},
is defined combinatorially (so we shall refer to it as Wan's
Hodge polygon in the remark below). We shall compare it
with classical Hodge polygons in the literature.

\begin{remark}
(i) Wan's Hodge polygon does not generally acquire
a geometric meaning, and we do not know of one for
the one-variable exponential sum case studied in this paper.
Nevertheless there is a well-known case in which it
does. If $f$ is an $n$-variable Laurent polynomial over a finite field,
Wan's Hodge polygon of the exponential sum of $f$
is defined in \cite[Section 1]{Wan:2}.
It is known that
if $f=0$ defines a toric variety (denoted by $X$)
then Wan's Hodge polygon of the exponential sum of $f$
does coincide with a variant of classical Hodge polygon which
is defined by `Hodge numbers' of
certain subgroup of the cohomology $H_c^{n-1}(X,\C)$
with compact support over the complex $\C$
(see \cite[Section 5]{Adolphson-Sperber:1}
for details and proofs).
This explains Wan's terminology of `Hodge polygon'.

(ii)
Let $X$ be a smooth projective
scheme of finite type over a finite field $\F_q$
such that the Hodge cohomology groups
$H^j(\hat{X},\Omega^i_{\hat{X}/W(\F_q)})$
are free $W(\F_q)$-modules of finite ranks $h^{i,j}(\hat{X})$,
where $W(\F_q)$ is the ring of Witt vectors over $\F_q$
and $\hat{X}/W(\F_q)$ is a lift of $X/\F_q$.
Recall that the ($m$-dimensional) classical Hodge polygon of $X$
(\`a la Katz-Mazur)
consists of line segments
of slope $i$ of horizontal length $h^i:=h^{i,m-i}(\hat{X})$
for all $i=0,\ldots,m$. For example,
the ($1$-dimensional) classical Hodge polygon of
a curve $X$ over a finite field of genus $g$ consists of a slope-$0$ segment
of length $h^{0}=g$ and a slope-$1$ one of length $h^{1}=g$.
It is known that this classical Hodge polygon of $X$
is a lower bound of its Newton polygon (see \cite{Katz,Mazur}).
However, for Artin-Schreier curves
this lower bound is not sharp and there is a sharp lower bound, that is
precisely Wan's Hodge polygon (after blown up by a factor of $p-1$
horizontally and vertically). In summary, Wan's Hodge polygon
is an analog of Katz-Mazur's classical Hodge polygon as
a lower bound to Newton polygons.

(iii)
From a geometrical point of view, it has long been a myth why
the Newton polygons of Artin-Schreier curves $X_f$
are much higher than the classical Hodge polygon
(see \cite{van der Geer} for example).
It is discovered recently
that for the exponential sum of
a one-variable {\em rational} function $f$
there is a generalized Wan's Hodge polygon
that determined by the orders of poles of $f$.
(This was conjectured by Poonen and Adolphson-Sperber
independently,
proved by \cite[Theorem 1.1]{Zhu:2}.)
Our Theorem \ref{T:3} below is asserting that
if $f$ is a polynomial then
Wan's Hodge polygon is asymptotically
(for $p$ large) a best lower bound!
We anticipate that a natural generalization of this sort
should hold true also for rational functions.
\end{remark}

\begin{remark}
The choice of the primitive $p$-th root of unity $\zeta_p$ does not
affect the Newton polygon. In fact, if $\zeta_p$ is replaced by
$\zeta_p^i$ in (\ref{E:S_ell}), then every coefficient $b_n$
in the $L$ function in (\ref{E:20})
will be changed by replacing every
$\zeta_p$ in its expression by $\zeta_p^i$.
The $p$-adic valuation
of $b_n$ is invariant under Galois conjugation.
\end{remark}

In this paper we prove Theorems  \ref{T:3} and \ref{T:NP}.
Part of Theorem \ref{T:3} was formulated as a conjecture by Daqing Wan
communicated to me in 2001 and it is now a one-dimensional case of a
new conjecture collected in \cite[Section 1.4]{Wan:2}.  The case $d=3$
of \ref{T:3} is proved by \cite[(3.14)]{Sperber:2} using Dwork's
method. Theorem \ref{T:3} also yields a complete answer
to a question (in one-variable case)
proposed by Katz on page 151 \cite[Chapter 5.1]{Katz:2}.
The first slope case is proved recently by \cite{SZ:2} by a
slope estimate technique essentially following Katz \cite{Katz:1}. A
weaker version of this theorem, which restricts to $f\in\cU(\Q)$, is proved
in \cite{Zhu:1} recently.
Recall from \cite[Section 5]{Zhu:1} that
 $\GNP(\A^d;\F_p):=\inf_{\bar{f}\in\A^d(\F_p)}\NP(\bar{f})$
if exists.

\begin{theorem}\label{T:3}
There is a Zariski dense open subset $\cU$ defined over $\Q$
in $\A^d$ such that
if $f\in\cU(\bar\Q)$ and if $\cP$ is a prime ideal in the ring of integers
of $\Q(f)$ lying over $p$, we have
for $p$ large enough (depending only on $f$),
$\NP(f\bmod \cP)=\GNP(\A^d;\F_p)$.
In particular, for every $f\in \cU(\bar\Q)$
one has $\lim_{p\rightarrow\infty}\NP(f\bmod\cP)=\HP(\A^d).$
\end{theorem}

In the proof of Theorem \ref{T:3} (Section \ref{S:proof}) we give
an explicit formula for
the asymptotic generic Newton polygon $\GNP(\A^d;\F_p)$
for $p$ large enough, which depends only on $d$ and the residue class
of $p\bmod d$. We consider Theorem
\ref{T:NP} as a major technical breakthrough of the present paper.
Since it is more involved we postpone
its discussion to Section \ref{S:NP}.
Our theorem has the following application in approximating
slopes of Artin-Schreier curves (see \cite[Corollary 1.3]{Zhu:1} for
a proof).

\begin{corollary}\label{C:1}
There exists a Zariski dense open subset $\cU$
defined over $\Q$ in $\A^d$ such that if
$f\in\cU(\bar\Q)$ and $\cP$ is any prime ideal in the ring
of integers of $\Q(f)$ lying over $p$, we have
$$
\lim_{p\rightarrow\infty}\frac{\NP(X_f\bmod \cP)}{p-1}
=\HP(\A^d).
$$
\end{corollary}

The paper is organized as follows. Section \ref{S:Dwork} carries on our
exposition of Dwork $p$-adic analysis from \cite[Section 2]{Zhu:1},
with emphasis on the semilinear theory, that is, the spot light is
at exponential sums of finite
fields that are {\em not} prime fields (in contrast to that in
\cite[Section 2]{Zhu:1}).
Section \ref{S:NP} contains the main technical theorem
(in Theorem \ref{T:NP}) of this paper,
in which we prove that for a large class
of matrix $F$ representing $\tau^{-1}$-linear Frobenius map
over a p-adic ring $\cO_a$ the  Newton polygon
of its characteristic polynomial coincides with that
of the matrix $F_a$ representing linear Frobenius map.
Section \ref{S:Zariski} defines
appropriate Zariski dense open subset
$\cW_r$ in $\A^{d-1}$ for each residue class $r \bmod d$.
Proofs for Theorems \ref{T:3}
lie in Section \ref{S:proof}. In the same section
we also prove that a certain stronger version of
Theorem \ref{T:3} is false, which answers a question of Daqing Wan.

\begin{acknowledgments}
  I am deeply indebt to Hanfeng Li and Daqing Wan
  for enlightening suggestions and warm encouragements.
  I also thank Alan Adolphson and
  the referee for comments on an earlier version.
  This research was partially
  supported by a grant of Bjorn Poonen from the David and Lucile
  Packard Foundation and
  the University of California at Berkeley.
\end{acknowledgments}

\section{Dwork $p$-adic theory in a nutshell}\label{S:Dwork}

The present section is in a sequel to \cite[Section 2]{Zhu:1}
yet it is self-contained for the convenience of the reader.
We formulates the Dwork trace formula following \cite{Dwork:1, Dwork:2,
Sperber:1,Sperber:2} at various stages without further notice.
Our trace formula of this article concerns
a Frobenius action on a
finite dimensional quotient space,
while that of \cite{Zhu:1}
is considering a Frobenius action on an infinite dimensional
vector space.

Let $\Q_{p^a}$ denote the unramified extension of
$\Q_p$ of degree $a$.
Let $\Omega_1=\Q_p(\zeta_p)$ and let
$\Omega_a=\Q_{p^a}(\zeta_p)$.
So $\Omega_a$ is the unramified extension of $\Omega_1$ of degree $a$.
Let $\cO_1=\Z_p[\zeta_p]$ and $\cO_a$ be the rings of integers
in $\Omega_1$ and $\Omega_a$, respectively.
Let $\tau$ be the lift of Frobenius endomorphism $c\mapsto c^p$
of $\F_q$ to $\Omega_a$ which fixes $\Omega_1$.

Fix $f(x)=x^d+\sum_{i=1}^{d-1}a_ix^i\in \bar{\Q}[x]$.
Suppose
$\bar{f}(x)=x^d+\sum_{i=1}^{d-1}\bar{a}_ix^i\in \F_q[x]$ is a
reduction of $f(x)$.
Let $\hat{f}(x)=x^d+\sum_{i=1}^{d-1}\hat{a_i}x^i$,
where $\hat{a}_i$ is the Teichm\"uller lifting of
$\bar{a}_i$, that is, $\hat{a}_i$ lies in $\Z_{p^a}$
such that
$\hat{a}_i\equiv \bar{a}_i\bmod p$ and $\hat{a}_i^q=\hat{a}_i$.
We shall write
$\vec{\bar{a}}:=(\bar{a}_1,\ldots,\bar{a}_{d-1})$ and
$\vec{\hat{a}}:=(\hat a_1,\ldots,\hat a_{d-1})$.

Let $\theta(x)=E(\gamma x)$ where $E(\cdot)$ is the $p$-adic
Artin-Hasse
exponential function and $\gamma$ is a root of $\log E(x)$ with
$\ord_p \gamma = \frac{1}{p-1}$ (as defined in Section 1).
We may write $\theta(x)=\sum_{m=0}^\infty \lambda_mx^m$
for $\lambda_m\in\cO_1$ . In fact, $\lambda_m=\frac{\gamma^m}{m!}$
and $\ord_p\lambda_m=\frac{m}{p-1}$ for $0\leq m\leq p-1$, and
$\ord_p\lambda_m\geq \frac{m}{p-1}$ for $m\geq p$.
Let $\vec{m}=(m_1,\ldots,m_{d-1})\in\Z_{\geq 0}^{d-1}$
and let $\vec{A}=(A_1,\ldots,A_{d-1})$ be variables.
Denote $\vec{A}^{\vec{m}}:=A_1^{m_1}\cdots A_{d-1}^{m_{d-1}}$.
Define for any $n\geq 0$ a polynomial in $\Omega_1[\vec{A}]$ below
$$G_n(\vec{A}):= \sum_{\substack{m_\ell\geq 0\\\sum_{\ell=1}^{d}\ell m_\ell
= n}}\lambda_{m_1}\cdots\lambda_{m_d}\vec{A}^{\vec{m}}.
$$
Let $G(X):=\prod_{i=1}^{d}\theta(\hat{a}_i X^i).$ So
$G(X)\in\cO_a[[X]]$ and its expansion is precisely
$G(X)=\sum_{n=0}^{\infty}G_n(\vec{\hat{a}})X^n\in\cO_a[[X]]$.

Let $K$ be a $p$-adic field over $\Omega_a$.
For any $c>0$ and $b\in\R$ let $\cL_K(c,b)$ be the set of power series
$\sum_{n=0}^{\infty}B_nX^n\in\Omega_a[[X]]$ with
$B_n\in K$ and $\ord_pB_n\geq cn+b$. Let
$\cL_K(c)=\bigcup_{b\in\R}\cL(c,b)$. For example, one may check
$G(X)\in\cL_K(\frac{1}{d(p-1)})$.
Note that $\cL_K(c)$ is a infinite dimensional vector space over $K$.

Consider the composition
$\alpha:=\tau^{-1}\cdot\psi\cdot G(X)$ on $\cL_K(c)$, where
$\psi$ is the Dwork $\psi$-operator on $\cL_K(c)$ defined by
$\psi(\sum_{n=0}^{\infty}B_nX^n)=\sum_{n=0}^{\infty}B_{np}X^n$,
and $G(X)$ denotes the multiplication map by $G(X)$.
One observes that
$\alpha$ is $\tau^{-1}$-linear $\Omega_a$-endomorphism of
$\cL_{\Omega_a}(\frac{p}{d(p-1)})$.
Write $\cL$ and $\cL^1$ for
$\cL_{\Omega_a}(\frac{p}{d(p-1)})$
and its subspace with no constant terms, respectively.
For $\ell\geq 0$ let
$\gamma_{\ell}:=\sum_{j=0}^{\ell}\frac{\gamma^{p^j}}{p^j}$.
Let
$$
R(X):=\sum_{\ell=0}^{\infty}\gamma_\ell \hat{f}^{\tau^\ell}(X^{p^\ell})
=
\sum_{\ell=0}^{\infty}
     \gamma_\ell\sum_{i=1}^{d}{\hat a_i}^{p^\ell}X^{ip^\ell}.
$$
Let $\nabla$ be a differential operator on $\cL$
defined formally by
$$\nabla:=\exp(-R(X))\cdot X\frac{\partial}{\partial X}\cdot\exp(R(X)).$$
For any $\sum_{n=0}^{\infty}B_nX^n\in\cL$,
we have
$$
\nabla(\sum_{n=1}^{\infty}B_nX^n) =\sum_{n=0}^{\infty}nB_nX^n
+(X\frac{\partial R(X)}{\partial
  X})(\sum_{n=0}^{\infty}B_nX^n).
$$
Clearly $\nabla(\cL)\subseteq \cL^1$ and we define
$\cM:=\cL^1/\nabla(\cL).$ Then $\cM$ has the induced
$\tau^{-1}$-linear endomorphism $\alpha$.
(See \cite[page 279]{Sperber:2} for more details.)

Let $\vec{e}$ denote the set of images of $\{X,X^2,\ldots,X^{d-1}\}$
in the quotient space $\cM$.
Then $\vec{e}$ form a basis for $\cM$
over $\Omega_a$, and
$\dim_{\Omega_a}\cM=d-1$.  Let $F$
be the matrix representation of $\alpha$ on $\cM$
with respect to the basis $\vec{e}$.
Let $G^{[a]}(X):= \prod_{j=0}^{a-1}G^{\tau^{j}}(X^{p^j}).$
Let $\alpha_a:=\psi^a\cdot G^{[a]}(X)$, which is a
(linear!) endomorphism of $\cM$ over $\Omega_a$.
The case $a=1$ is thoroughly studied
in \cite{Zhu:1}. Let $F_a$ be the matrix
representation of $\alpha_a$ on $\cM$ with respect to
this monomial basis $\vec{e}$.
The map $\alpha$ on $\cM$ is given by
$\alpha\vec{e}=\vec{e}F$. Since $\alpha_a=\alpha^a$ and
$\alpha$ is $\tau^{-1}$-linear we see easily
that
$F_a=FF^{\tau^{-1}}\cdots F^{\tau^{-(a-1)}}.$

For any positive integer $n$ let $\M_n(\cdot)$ denote the set of
all $n$ by $n$ matrix over some ring.  Let $\I_n$ denote the $n$
by $n$ identity matrix. By Dwork trace formula (see \cite[Theorem
2.2]{Hong:1} or \cite[Section 2 and in particular
(2.35)]{Sperber:2} \cite[discussions in Section 2]{Dwork:1} for
details) for any prime $\cP$ over $p$ of degree $a$ we have
\begin{eqnarray}\label{E:Dwork}
L(f\bmod\cP;T)=\det(\I_{d-1}-T\alpha_a|\cM)
=\det(\I_{d-1}-TF_a).
\end{eqnarray}

\begin{remark}\label{R:diff}
One observes that the computation of the above $L$ function is
reduced to the process of diagonalization (or triangularization) of
the matrix $F_a$. Write $\Q_{p^\infty}$ for the fraction field of
$W(\bar\F_p)$.
Even though the Dieudonn\'e-Manin classification \cite{Manin:1}
asserts that it is plausible over
$\Q_{p^\infty}(\zeta_p)$,
it is far more than a small business in practice.
One of our hardest tasks broils down to proving
a stronger version of the
Dieudonn\'e-Manin classification holds in the sense that
our matrix can be diagonalized  (or triangularized) over
the base field $\Omega_a$. This is accomplished in Section \ref{S:NP}.
Our other challenging tasks include finding the Zariski
open subset set defined over $\Q$, which is done in
Section \ref{S:Zariski}.
\end{remark}

\begin{example}
Below we give a simple example only to demonstrate the essential
difficulty and new effects amounted in Wan's conjecture
when passing from $\F_p$ to $\F_{p^a}$.
Let $p\equiv -1\bmod 4$ and $f(x)=x^4+cx$ in $\F_{p^a}$.
For $p$ large enough we can compute and get
$$F =
\left(
\begin{array}{ccc}
  \frac{\gamma^{p-1}}{(p-1)!}\hat{c}^{p-1}
& \frac{\gamma^{p-2}}{(p-2)!}\hat{c}^{p-2}
& \frac{\gamma^{p-3}}{(p-3)!}\hat{c}^{p-3}
\\
  \frac{\gamma^{2p-1}}{(2p-1)!}\hat{c}^{2p-1}
& \frac{\gamma^{2p-2}}{(2p-2)!}\hat{c}^{2p-2}
& \frac{\gamma^{2p-3}}{(2p-3)!}\hat{c}^{2p-3}
\\
\frac{\gamma^{3p-1}}{(3p-1)!}\hat{c}^{3p-1}
&\frac{\gamma^{3p-2}}{(3p-2)!}\hat{c}^{3p-2}
&\frac{\gamma^{3p-3}}{(3p-3)!}\hat{c}^{3p-3}
\end{array}
\right)$$
in $\M_{3}(\Omega_a)$,
where $\hat{c}$ is the Teichm\"uller lifting of $c$ in
$\Omega_a$. Then our
Frobenius matrix is $F_a =FF^{\tau^{-1}}\cdots
F^{\tau^{-(a-1)}}$.
If $a=1$, the diagonalization process is linear algebra,
that is, one needs $C^{-1}FC$ diagonal for some $C$
over $\Omega_1$.
For $a>1$, the process is $\tau^{-1}$-linear,
that is, one needs $C^{-\tau}FC$ diagonal for some $C$
over $\Omega_a$.
The semilinear algebra involved is highly nontrivial
(see Proposition \ref{P:keylemma}).
The reader who is interested in complete
numerical analysis of lower degree cases
are referred to two new papers
\cite{Hong:1} and \cite{Hong:2}.
\end{example}

\section{The two Newton polygons}\label{S:NP}

This section is technical and a key technical ingredient
in our argument is a version of
$p$-adic Banach fix point theorem.

Let $m$ be a positive integer.  For any $m$ by $m$ matrix $M$ in
with coefficients in $\Omega_a$ (i.e., $M\in \M_{m}(\Omega_a)$)
and any $1\leq n\leq m$
let $M^{[n]}$ denote the submatrix of $M$
consisting of its first $n$ rows and columns.
Let $M_a:=MM^{\tau^{-1}}\cdots M^{\tau^{-(a-1)}}.$
We observe that if
$\ord_pM_{ij} \rightarrow \ord_pM_{i1}$ for every $j$,
$\ord_pM_{i,1}-\ord_pM_{i-1,1}>\xi$ for some constant $\xi>0$ for every $i$,
and $\ord_p\det M^{[n]}\rightarrow \sum_{i=1}^{n}\ord_pM_{i1}$
for every $1\leq n\leq m$, then
$\NP_q(\det(\I_m-TM_a))=\NP_p(\det(\I_m-TM))$.
This observation is highly nontrivial, so we will
prove it in \ref{P:keylemma} below.
Let
\begin{eqnarray*}
\delta(M)&:=&(p-1)\min_{1\leq i\leq m-1} (\min_{1\leq j\leq
m}\ord_pM_{i+1,j}-
 \max_{1\leq j\leq m}\ord_pM_{ij});\\
\eta(M)&:=&(p-1)\max_{1\leq n\leq m-1}(\ord_p\det M^{[n]}
                -\sum_{i=1}^{n}\min_{1\leq j\leq n+1}\ord_pM_{ij}).
\end{eqnarray*}
It is easy to observe that $\eta$ and $\delta$ are
nonnegative integers.

Write $\det(\I_m-TM)=1+c'_1T+\cdots+c'_mT^m\in \Omega_a[T]$,
and $\det(\I_m-TM_a)=1+b'_1T+\cdots+b'_mT^m$.

\begin{proposition}\label{P:keylemma}
Let $M$ be in $\M_m(\cO_a)$ (recall $\cO_a$ is the ring of integers in
$\Omega_a$) such that $\delta(M)> m\,\eta(M)$.
Then $\ord_pc'_n=\ord_p\det M^{[n]}.$
There exists a unique upper triangular matrix $C$ in $\M_m(\Omega_a)$ with all
$1$'s on its diagonal and with $\ord_p C_{ij}\geq -\frac{\eta}{p-1}$
such that $M':=C^{-\tau}MC$ in $\M_m(\Omega_a)$
is lower triangular.  Set $\det M^{[0]}: = 1$.
For any $1\leq n\leq m$ one has
\begin{eqnarray}\label{E:M'}
\ord_pM'_{nn}&=&\ord_p\det M^{[n]}-\ord_p\det M^{[n-1]}.
\end{eqnarray}
Moreover, $M'$ has
strictly increasing $p$-adic orders down its diagonal.
\end{proposition}

To be useful to the reader, we make some remarks on what
leads us to the formulation of the hypothesis in
\ref{P:keylemma} and \ref{T:NP}: Choose a different basis
$\vec{e}_w:=\{(\gamma^{1/d}X)^i\}_{1\leq i\leq d-1}$
for
$\cL_{\Omega_a(\gamma^{1/d})}(\frac{p}{d(p-1)})$
over $\Omega_a(\gamma^{1/d})$.
Let $F^w$ be the matrix for $\alpha$  under this basis.
Then one notes that
$\ord_p F^w_{ij}\geq i/d+ r_{ij}/d(p-1)$ where
$r_{ij}$ is the least nonnegative residue of
$-(pi-j)\bmod d$.
As $p\rightarrow\infty$
this lower bound of $\ord_pF^w_{ij}$ converges to $i/d$ for
every $i$ and $j$.

\begin{proof}
1)
Let $C$ be an upper triangular matrix with
$(i,j)$-th entry denoted by indeterminant $C_{ij}$
and with all $1$'s on its diagonal.
For any $1\leq j\leq m$ and
all $i=1,\ldots, j-1$ set
\begin{eqnarray}\label{E:1}
(C^{-\tau}MC)_{ij}&=&0.
\end{eqnarray}
Write $D:=C^{-\tau}$. Then $D$ is also upper
triangular with all $1$'s on its diagonal.
So we have
\begin{equation*}
(C^{-\tau}MC)_{ij}
=\sum_{k=i}^{m}\sum_{\ell=1}^{j}D_{ik}M_{k\ell}C_{\ell j}
=\sum_{\ell=1}^{j-1}M_{i\ell}C_{\ell j} + M_{ij}
 +\sum_{k=i+1}^{m}\sum_{\ell=1}^{j}M_{k\ell}(D_{ik}C_{\ell j})=0.
\end{equation*}
Then one verifies that
(\ref{E:1})
is equivalent to
\begin{equation}\label{E:6}
M^{[j-1]}\left(
\begin{array}{c}
C_{1j}\\
\vdots\\
C_{j-1,j}
\end{array}
\right)
+
\left(
\begin{array}{c}
M_{1j}\\
\vdots\\
M_{j-1,j}
\end{array}
\right)
+
\left(
\begin{array}{c}
\sum_{k=2}^{m}\sum_{\ell=1}^{j}M_{k\ell}(D_{1k}C_{\ell j})\\
\vdots\\
\sum_{k=j}^{m}\sum_{\ell=1}^{j}M_{k\ell}(D_{j-1,k}C_{\ell j})
\end{array}
\right)
=0.
\end{equation}

Now we introduce some notations.
For $1\leq i,j\leq m$,
let $M_{(i,j)}$ denote the submatrix of $M$ with its $i$-th row and
the $j$-th column removed.
Let $M^*$ denote the {\em adjoint matrix} of $M$, that is,
the matrix whose $(i,j)$-th entry is equal to $(-1)^{i+j}\det
M_{(j,i)}$.  {}From linear algebra
we have that $MM^*=M^*M=\det M$.
Consider $\vec{C}$-monomials,
i.e., consider all $C_{ij}$'s and $C_{ij}^\tau$'s as variables
where $1\leq i<j\leq m$.
Our hypothesis on $M$ implies
that $\det M^{[n]}\neq 0$ for every $1\leq n\leq m$.
Thus we may multiply $(M^{[j-1]})^{-1}$ on
the left-hand-side of (\ref{E:6}) and get
for all $1\leq i<j$
\begin{eqnarray}
\label{E:4}
C_{ij} &=& w_i(\vec{C})+ v_i
\end{eqnarray}
where
$w_i(\vec{C})=-\frac{(M^{[j-1]})^*}{\det(M^{[j-1]})}
\sum_{k=i+1}^{m}\sum_{\ell=1}^{j}M_{k\ell}(D_{ik}C_{\ell j})$
and
$v_i=-\frac{(M^{[j-1]})^*}{\det(M^{[j-1]})}M_{ij}$.
It is easy to see that
\begin{eqnarray}\label{E:vi}
\ord_p (M^{[j-1]})^* M_{ij}
&\geq &\min_{1\leq i\leq j-1}\ord_p\sum_{\ell=1}^{j-1}
(-1)^{i+\ell}
(\det M^{[j-1]}_{(\ell,i)}) M_{\ell j}\\\nonumber
& \geq &
\sum_{\ell=1}^{j-1}\min_{1\leq k\leq j}
(\ord_pM_{\ell k}).
\end{eqnarray}
Thus
\begin{equation}\label{E:latter}
\ord_p(v_i)
\geq
 \sum_{\ell=1}^{j-1}\min_{1\leq k\leq j}
(\ord_pM_{\ell k})-\ord_p\det M^{[j-1]}
\geq -\frac{\eta}{p-1}.
\end{equation}
Then the $p$-adic valuation of
coefficients of any $\vec{C}$-monomial in
$\vec{w}(\vec{C})$
is $\geq\frac{\delta-\eta}{p-1}$ by comparing to $\vec{v}$.

Now change variables by
setting $X_{ij}:=\gamma^{\eta}C_{ij}$
for all $1\leq i<j\leq m$.
Write $z_i(\vec{X})$ for $\gamma^\eta w_i(\vec{C})$ as polynomials
in $X_{ij}$ and $X_{ij}^\tau$,  one has
$X_{ij} = z_i(\vec{X})+v_i\gamma^\eta$.
We claim that the right-hand-side of
(\ref{E:4}) for all $j=1,\ldots, m$ together
defines a contraction map
with regard to $X_{ij}$'s
on $\cO_a^{\frac{m(m-1)}{2}}$.
It suffices to show that $z_i(\vec{X})$ has all coefficients
of $p$-adic valuation positive. If this is the case,
then the Banach fixed point theorem applies and
one  has integral solutions $X_{ij}$
and consequently $C,M'\in\M_m(\Omega_a)$
with $\ord_pC_{ij}\geq -\frac{\eta}{p-1}.$

Note that $D=C^{-\tau}=(C^*)^{\tau}$.
For any $1\leq t<k\leq m$, it is an exercise to show that
$D_{tk}=(-1)^{t+k}\det(C_{(k,t)})^\tau$
is a degree $k-t$ polynomial in $\vec{C}$.
Thus $w_i(\vec{C})$ is of degree $\leq m$ in $\vec{C}$
Now it is another elementary exercise to
show that coefficients of
$\vec{X}$ in $z_i(\vec{X})$
has $p$-adic order
$\geq \frac{\delta - m\eta}{p-1}$,
which is positive by our hypothesis upon $M$.

2) From now on we assume $C$ is as chosen above.
It is an exercise to see for all $1\leq j\leq m$ one has
$$
M'_{jj}
 = (DMC)_{jj}
 =\sum_{1\leq i<j\leq m}M_{ji}C_{ij}+M_{jj}+
\sum_{1\leq i\leq j<k\leq m}M_{ki}D_{jk}C_{ij}.
$$
Also note that
\begin{equation*}
\det M^{[j]}= \det M^{[j-1]}
(M_{jj}+\sum_{i=1}^{j-1}M_{ji}v_i),
\end{equation*}
so one has
\begin{equation}\label{E:13}
M'_{jj}= \frac{\det M^{[j]}}{\det M^{[j-1]}}
+\sum_{i=1}^{j-1}M_{ji}w_i
+\sum_{1\leq i\leq j<k\leq m}M_{ki}D_{jk}C_{ij}.
\end{equation}
By some simple computations,
one finds every term exact the first one
on the right-hand-side of (\ref{E:13})
has
$p$-adic valuation $>\ord_p\det M^{[j]}-\ord_p \det M^{[j-1]}$.
Applying the isoscele
principle, one concludes that
$\ord_p(M'_{jj})= \det M^{[j]}-\ord_p\det M^{[j-1]}.$
This proves (\ref{E:M'}).

3)
By a similar argument as above,
one can show that there exists an upper triangular matrix
$C'$ with all $1$'s on the diagonal such that
$M'':={C'}^{-1}MC'$ is lower triangular and
\begin{eqnarray*}
\ord_pM''_{jj} &=&
\ord_p \det M^{[j]} - \ord_p \det M^{[j-1]}
\end{eqnarray*}
for $1\leq j\leq m$.
It follows easily that
$\ord_p c'_n = \sum_{\ell=1}^{n}\ord_p M''_{\ell\ell}
=\ord_p\det M^{[n]}$.

4)
Finally we shall omit the proof of the last statement.
The basic idea is using (\ref{E:M'}) to reduce to
show
$
\ord_p \det M^{[n+1]} + \ord_p \det M^{[n-1]} >
2\,\ord_p \det M^{[n]}.
$
\end{proof}

\begin{remark}\label{R:2-polygons}
For the purpose of computing $L$ function
according to Dwork's trace formula (\ref{E:Dwork}),
the relation between $\NP_q(\det(\I_{d-1}-TF_a))$ and
$\NP_p(\det(\I_{d-1}-TF))$ has been
explored in the literature (see for example \cite{Wan:1})
since the latter is much more straightforward to
compute. However, as passing from $F$ to $F_a$, the only thing
we knew previously is that their corresponding Newton polygons
have the same lower bounds (i.e., the Hodge polygon)
and these Newton polygons are not
generally equal.
See some discussion including a good example
in \cite[Section 1.3]{Katz}
and a study of ordinary case in \cite[Theorem 2.4]{Wan:1}. Little
is known besides these, yet the passage of Newton polygon data
from $F$ to $F_a$ is the bottleneck in sharp slope estimations
generally. In the theorem below we formulate an explicit
criterion under which the two aforementioned Newton polygons
coincide.
\end{remark}

\begin{theorem}\label{T:NP}
Let $M$ be in $\M_m(\cO_a)$ such that $\delta(M)>m\,\eta(M)$.
Then $\NP_q(\det(\I_m-TM_a))=\NP_p(\det(\I_m-TM))$,
and they are
equal to the lower convex hull in $\R^2$ of the
points $(n,\ord_p\det M^{[n]})$ for $0\leq n\leq m$.
\end{theorem}
\begin{proof}
Let $C, M' \in\M_{d-1}(\Omega_a)$ be as in \ref{P:keylemma},
then $C^{\tau^a}=C$ and ${M'}^{\tau^a} = M'$.
So
\begin{eqnarray*}
C^{-\tau}M_aC^{\tau}
&=& C^{-\tau^{a+1}}
    (M^{\tau^a} M^{\tau^{a-1}}\cdots M^{\tau})C^{\tau}\\
&=& (C^{-\tau} M C)^{\tau^{a}}(C^{-\tau} M C)^{\tau^{a-1}}
    \cdots (C^{-\tau} M C)^{\tau}\\
&=& {M'}^{\tau^{a}} {M'}^{\tau^{a-1}}\cdots {M'}^{\tau}\\
&=& (M')_a.
\end{eqnarray*}
By \ref{P:keylemma}, one knows that $(M')_a$ is lower triangular
with
\begin{equation*}
\ord_p((M')_a)_{jj}
= \sum_{\ell=0}^{a-1}\ord_p (M'_{jj})^{\tau^\ell}
= a(\ord_p M'_{jj}) = a (\ord_p\det M^{[j]} - \ord_p\det M^{[j-1]}).
\end{equation*}
Since
$$\det(\I_m-T M_a)
=\det(\I_m - T C^{-\tau} M_a C^{\tau})
=\det(\I_m-T (M')_a),
$$
and that $\ord_pM'_{jj}$ being strictly increasing according to $j$,
one has
\begin{eqnarray*}
\ord_qb'_n
&=&\frac{1}{a}\sum_{j=1}^{n}\ord_p((M')_a)_{jj}
=\frac{1}{a}\sum_{j=1}^{n}a(\ord_p\det M^{[j]} - \ord_p\det M^{[j-1]})\\
&=&\ord_p(\det M^{[n]}).
\end{eqnarray*}
Thus $\ord_qb'_n
=\ord_p\det M^{[n]}
=\ord_q c'_n$
for every $n$.
This finishes the proof.
\end{proof}

\section{Zariski dense open subset $\cW_r$ in $\A^{d-1}$}
\label{S:Zariski}

We shall use an auxiliary $d-1$ by $d-1$
matrix $F^\dagger$ defined by
$F^\dagger_{ij}(\vec{\hat a}):=G^{\tau^{-1}}_{pi-j}(\vec{\hat a})$
for every $\vec{\hat a}$. We outline our approach as below:
(1) we find a Zariski dense open subset $\cX_r$ of $f$'s in which
$\ord_pF_{ij} = \ord_pF_{ij}^\dagger =
\frac{\lceil\frac{pi-j}{d} \rceil}{p-1}$;
(2) we find a Zariski dense open subset $\cW_r$ of $f$'s in which
$\ord_p(\det F^{[n]})= \ord_p(\det(F^\dagger)^{[n]})
=\frac{n(n+1)}{2d}+\epsilon_n$,
both for all $p\equiv r\bmod d$ and $p$ large enough.
Basically we are looking for sufficient condition on
$p$ and $f$ such that the Frobenius matrix $F$
satisfies the hypothesis of \ref{P:keylemma}.
This is the key observation
prepared for the proof in Section \ref{S:proof}.

We adopt the same notation as that in \cite[Section 3]{Zhu:1}.
For convenience of the reader, we give complete definitions
for all statements of our theorems.
Let $r$ be a positive integer with $1\leq r\leq d-1$ and
$\gcd(r,d)=1$ for the rest of the section.
For any $1\leq i,j\leq d-1$, let $r_{ij}$
(resp.  $r'_{ij}$) be the least nonnegative residue of
$-(ri-j) \bmod d$ (resp. $ri-j\bmod d$).
Let $\delta_{ij}=0$ for $j<r'_{i1}+1$ and let $\delta_{ij}=1$ for
$j\geq r'_{i1}+1$.

Let $1\leq n\leq d-1$.
Let $\vec{v}:=(v_1,\ldots,v_n)\in\Z_{\geq 0}^n$
let $|\vec{v}|:=\sum_{\ell=1}^{n}v_\ell$
and $\vec{v}!:=v_1!\cdots v_{n}!$.
For any $0\leq t\leq n$
let $S_n^t$ denote the subset of the
symmetric group $S_n$ consisting of all
$\sigma$ such that
$\sum_{i=1}^{n}r_{i,\sigma(i)}
=\min_{\sigma'\in S_n}\sum_{i=1}^{n}r_{i,\sigma'(i)} + dt$.
For any $1\leq i,j\leq n$ and
$0\leq s\leq n$ define a subset of $\Z_{\geq
0}^{d-1}$ by
\begin{eqnarray*}
\cM_{ij}^s:=
\{\vec{m}=(m_1,\ldots,m_{d-1})\in\Z_{\geq 0}^{d-1}|
\sum_{\ell=1}^{d-1}\ell m_{d-\ell} = r_{ij}+ds
\}.
\end{eqnarray*}
Then let
$$
\begin{aligned}
H^s_{ij}(\vec{A})
&:=\sum_{\vec{m}\in\cM^s_{ij}}
\frac{(\frac{r_{i1}-1}{d}+d-1)(\frac{r_{i1}-1}{d}+d-2)\cdots
(\frac{r_{i1}-1}{d}-\delta_{ij}+s+1-|\vec{m}|)}{\vec{m}!}
\vec{A}^{\vec{m}}.
\end{aligned}
$$
Clearly $H^s_{ij}$ lies in $\Q[\vec{A}]=\Q[A_1,\ldots,A_{d-1}]$.
For any $0\leq t\leq n$ let
\begin{eqnarray*}
f_n^t(\vec{A})&:=&
\sum_{\substack{s_0+s_1+\cdots+s_n=t\\s_0,\ldots,s_n\geq 0}}
\sum_{\sigma\in S_n^{s_0}}\sgn(\sigma)
\prod_{i=1}^{n}H^{s_i}_{i,\sigma(i)}(\vec{A}).
\end{eqnarray*}
Let $t_n$ be the least nonnegative integer $t$ such that
$f_n^t\neq 0$.
Let
\begin{equation*}
\Psi_{d,r}(\vec{A}):=\prod_{0\leq j-1\leq i\leq d-1}H_{ij}^0(\vec{A}),
\quad\quad
\Phi_{d,r}(\vec{A}):=\prod_{1\leq n\leq d-1}f_n^{t_n}(\vec{A}).
\end{equation*}
Let $\cX_r$ and $\cY_r$ be the subset of $\A^{d-1}$ consisting of all
$f(x)=x^d+a_{d-1}x^{d-1}+\cdots+a_1x$ with
$\Psi_{d,r}|_{\vec{A}=\vec{a}}\neq 0$ and $\Phi_{d,r}|_{\vec{A}=\vec{a}}\neq
0$, respectively.  Let $\cW_r:=\cX_r\cap\cY_r$.

For any $b\in\Z$ let $\gamma^{>b}$ denote a term with
$\ord_p(\cdot)>\frac{b}{p-1}$. We define $\gamma^{\geq b}$ similarly.
For any $m$ by $m$ matrix $M$
and $n\leq m$ let $M^{[n]}$ denote the truncated submatrix of $M$
consisting of its first $n$ rows and columns.
Let
\begin{eqnarray}\label{E:epsilon_n}
\epsilon_n&:=&\frac{\min_{\sigma\in S_n}\sum_{\ell=1}^{n}
r_{\ell,\sigma(\ell)} + dt_n}
{d(p-1)}.
\end{eqnarray}

\begin{lemma}\label{L:Astar_Bound}
1) Let $\vec{a}\in \cX_r(\bar\Q)$.
There exists $N>0$ such that for $p>N$ we have for any
$0\leq j-1\leq i\leq d-1$ that
\begin{eqnarray}\label{E:G_Bound1}
\ord_pF^\dagger_{ij}(\vec{\hat a})
&=&\frac{\llceil\frac{pi-j}{d}\rrceil}{p-1}.
\end{eqnarray}

2) Let $\vec{a}\in \cY_r(\bar\Q)$.
There exists $N>0$ such that for $p>N$ we have
for every $1\leq n\leq d-1$ that
\begin{eqnarray}
\label{E:G_Bound2}
\ord_p\det (F^\dagger)^{[n]}(\vec{\hat a})
&=&\frac{n(n+1)}{2d}+\epsilon_n.
\end{eqnarray}
\end{lemma}
\begin{proof}
Since $F^\dagger_{ij}=G_{pi-j}^{\tau^{-1}}$, it suffices to
prove our assertion for $G_{pi-j}$.
Let notation be as in
\cite[4.2 and 4.3]{Zhu:1}.
By \cite[4.2]{Zhu:1},
for $p\geq (d^2+1)(d-1)$ we have
\begin{eqnarray}\label{E:H-K}
K^0_{ij}(\vec{A})&=&u_nH^0_{ij}(\vec{A})+\gamma^{\geq p-1}
\end{eqnarray}
for some $p$-adic unit $u_n$ in $\Z_p$.
Thus by \cite[4.3]{Zhu:1} one has
$$
G_{pi-j}(\vec{A})=u_nH^0_{ij}(\vec{A})
\gamma^{\llceil\frac{pi-j}{d}\rrceil}
+\gamma^{>\llceil\frac{pi-j}{d}\rrceil}.
$$
By the hypothesis $\vec{a}\in\cX_r(\bar\Q)$,
we have $H^0_{ij}(\vec{a})\neq 0$. So for
$p$ large enough one gets
\begin{eqnarray*}
\ord_p(G_{pi-j}(\vec{\hat a})) = \frac{\llceil\frac{pi-j}{d}\rrceil}{p-1}.
\end{eqnarray*}
This proves 1).
Part 2)  follows immediately from \cite[4.3]{Zhu:1}.
\end{proof}

\begin{lemma}\label{L:10}
Let $\vec{a}\in\cX_r(\bar\Q)$. For $0\leq j-1\leq i\leq d-1$ and for
$p$ large enough one has
\begin{eqnarray*}
\ord_p F_{ij}(\vec{\hat a})&=&\frac{\llceil\frac{pi-j}{d}\rrceil}{p-1}.
\end{eqnarray*}
\end{lemma}
\begin{proof}
The auxiliary matrix $F^\dagger$ is
$p$-adically close to $F$ in the following sense.
For any $\vec{a}\in\A^{d-1}(\bar\Q)$ and $1\leq i,j\leq d-1$ we have
\begin{eqnarray}\label{E:diff}
\ord_p(F_{ij}(\vec{\hat a})-F^\dagger_{ij}(\vec{\hat a}))
&\geq&\frac{pi-j}{d(p-1)}+\frac{p}{d(p-1)}.
\end{eqnarray}
(See \cite[Lemma 3.2]{Hong:1} for a complete proof
or follow the proof of Theorem 3.10 in \cite{Sperber:2}.)

By (\ref{E:G_Bound1}) and (\ref{E:diff})
there exists $N>0$ such that for all $p>N$ we have
\begin{eqnarray*}
\ord_p(F_{ij}(\vec{\hat a})-F^\dagger_{ij}(\vec{\hat a}))
&>&\ord_pF^\dagger_{ij}(\vec{\hat a}).
\end{eqnarray*}
By the isosceles triangle principle, we have
$\ord_p F_{ij}(\vec{\hat a})=\ord_p F^\dagger_{ij}(\vec{\hat a}),$
hence our assertion follows from (\ref{E:G_Bound1}).
\end{proof}

\begin{proposition}\label{P:10}
The subset $\cW_r$ is Zariski dense open in $\A^{d-1}$
defined over $\Q$.
For $\vec{a}\in\cW_r(\bar\Q)$ and $p$ large enough, one has
for all $1\leq n\leq d-1$  and $0\leq j-1\leq i\leq d-1$ that
\begin{equation*}
\ord_p F_{ij}(\vec{\hat a})=\frac{\llceil\frac{pi-j}{d}\rrceil}{p-1},
\quad \ord_p \det F^{[n]}(\vec{\hat a})=\frac{n(n+1)}{2d}+\epsilon_n.
\end{equation*}
\end{proposition}
\begin{proof}
By \cite[Section 3]{Zhu:1} one knows that
$\Psi_{d,r}\neq 0$ and $\Phi_{d,r}\neq 0$.
Thus the first assertion follows.
The first equality is precisely proved in Lemma \ref{L:10} above.
We shall focus on the second equality for the rest of our proof.
Write $\Delta$ for the set $\{1,\ldots,n\}$.
By definition,
\begin{eqnarray*}
\det F^{[n]}
&=& \sum_{\sigma\in S_n}\sgn(\sigma)
    \prod_{\ell=1}((F_{\ell,\sigma(\ell)} - F^\dagger_{\ell,\sigma(\ell)} )
                              +F^\dagger_{\ell,\sigma(\ell)}  )\\
&=& \sum_{\sigma\in S_n} \sgn(\sigma) \prod_{\ell=1}^{n}
                         F^\dagger_{\ell,\sigma(\ell)}\\
&&+ \sum_{\sigma\in S_n}\sgn(\sigma)\sum_{\Delta_1\subsetneq\Delta}
(\prod_{\ell\in\Delta_1}
(F_{\ell,\sigma(\ell)} - F^\dagger_{\ell,\sigma(\ell)})
\prod_{\ell'\in\Delta-\Delta_1} F^\dagger_{\ell',\sigma(\ell')}).
\end{eqnarray*}
By (\ref{E:diff}) and (\ref{E:G_Bound1})
(since $\vec{a}\in\cX_r(\bar\Q)$), for $p$ large enough we have
\begin{eqnarray*}
&&\ord_p(\det F^{[n]} - \det (F^\dagger)^{[n]})\\
&\geq&
\min_{\sigma\in S_n, \Delta_1\subsetneq \Delta}
(\sum_{\ell\in \Delta_1}
\ord_p(F_{\ell,\sigma(\ell)} - F^\dagger_{\ell,\sigma(\ell)})
+
\sum_{\ell'\in \Delta-\Delta_1}
\ord_pF^\dagger_{\ell',\sigma(\ell')}
) \\
&\geq&
\min_{\sigma\in S_n, \Delta_1\subsetneq \Delta}
(\sum_{\ell\in\Delta_1}
(\frac{p\ell-\sigma(\ell)}{d(p-1)}+\frac{p}{d(p-1)})
+\sum_{\ell'\in\Delta-\Delta_1}\frac{p\ell'-\sigma(\ell')}{d(p-1)})\\
&\geq&
\frac{n(n+1)}{2d}+\frac{p}{d(p-1)}.
\end{eqnarray*}
Since $\vec{a}\in\cY_r(\bar\Q)$ and
since $\epsilon_n$ goes to $0$ as $p$ approaches $\infty$,
for $p$ large enough this is strictly
greater than
$\ord_p\det (F^\dagger)^{[n]}$ by (\ref{E:G_Bound2}).
By the isosceles principle, we concludes our assertion.
\end{proof}

\section{The asymptotic generic Newton polygon
and slope filtration}\label{S:proof}

Let notations be as in previous sections.
In particular, recall the matrices $F$ and $F_a$ represent
the $\tau^{-1}$-linear and linear Frobenius maps
$\alpha$ and $\alpha_a$, respectively.

Let $\cW:=\bigcap_r\cW_r$ where $r$ ranges in
all $1\leq r\leq d-1$ with $\gcd(r,d)=1$.
Let $f(x)=x^d+a_{d-1}x^{d-1}+\cdots+a_1x+a_0\in
\A^d$.  Let $\cU$ be the pre-image of $\cW$ in $\A^d$ under the
projection map $\iota: \A^{d}\rightarrow\A^{d-1}$ by
$\iota(f)=\vec{a}$ with $\vec{a}=(a_1,\ldots,a_{d-1})$.
Let $\cP$ be any degree-$a$ prime ideal in $\Q(f)$ over $p$.
By (\ref{E:Dwork}) one has
\begin{eqnarray*}
L(x^d+\cdots +a_1x\bmod\cP;T)&=&\det(\I_{d-1}-F_a(\vec{\hat a})\cdot T)\\
&=&1+b_1(\vec{a})T+\cdots+b_{d-1}(\vec{a})T^{d-1}\in\Z[\zeta_p][T].
\end{eqnarray*}
One observes (see \cite[Section 5]{Zhu:1}) that
$\ord_q(b_n(f))=\ord_q(b_n(\vec{a}))$.

\begin{proof}[Proof of Theorem \ref{T:3}]
By \cite[Theorem 5.1]{Zhu:1}, for $p$ large enough one has
$\GNP(\A^d;\F_p)$
equal to the the lower convex hull of
points
\begin{eqnarray}\label{E:GNP}
(n,\frac{n(n+1)}{2d}+\epsilon_n)\quad \mbox{for $0\leq n\leq d-1$},
\end{eqnarray}
each of which is a vertex (recall $\epsilon_n$ from (\ref{E:epsilon_n})).
By \ref{P:10}, one sees that $\cW$ is Zariski dense open in $\A^{d-1}$ and
so is $\cU$ in $\A^{d}$ by its definition.
According to  the discussion preceding the proof,
it then suffices to prove our statements for $\vec{a}\in\cW(\bar\Q)$.
Suppose $\vec{a}\in\cW(\bar\Q)$.
{}From \ref{P:10} it is not hard to verify that
as $p$ increases $\delta(F)$ is unbounded
while $\eta(F)$ is bounded.
Thus for $p$ large enough,
$F(\vec{\hat a})\in\M_{d-1}(\Omega_a)$ clearly
satisfies the hypothesis of \ref{T:NP}.
So, by \ref{T:NP}, one has $\ord_qb_n(\vec{a})=\ord_p(\det F^{[n]})$.
Then by \ref{P:10} one sees that $\NP(f\bmod \cP)$ is equal to
the aforementioned convex hull given in (\ref{E:GNP}).
Finally, one notes that $\epsilon_n$ goes to $0$ as $p$
goes to infinity, this proves the theorem.
\end{proof}

\begin{remark}\label{R:NP}

We remark on another consequence of our results in
\ref{P:keylemma} and \ref{T:NP} from a different viewpoint.
It can be shown, by a symmetric argument
as that in \ref{P:keylemma},
that the two Frobenius matrices
$F$ and $F_a$ for  the exponential sums of $f$
are diagonalizable over the base field $\Omega_a$
provided $f\in\cU(\bar\Q)$ and $p$ is large enough.
This implies that for $f\in\cU(\bar\Q)$ and $p$ large enough,
the $F$-crystal $\cM$
(arisen from one-variable exponential sum in Section \ref{S:Dwork})
has a slope filtration over $\Omega_a$.
More precisely, it is isogenous over $\Omega_a$ to the direct sum of
rank-one $F$-crystals of slopes
$\frac{1}{d}+\epsilon_1$, $\frac{2}{d}+(\epsilon_2-\epsilon_1)$,
$\cdots$, and $\frac{d-1}{d}+(\epsilon_{d-1}-\epsilon_{d-2})$
respectively (necessarily in a strictly increasing order).
This provides an improvement, in
the case of one-variable exponential sums,
to the Dieudonn\'e-Manin classification
which asserts that the $F$-crystal has a slope filtration
over $\Q_{p^\infty}(\zeta_p)$
(see the classic of Manin
\cite[Chapter II]{Manin:1}, or see
\cite[Theorem 5.6]{Kedlaya:1} and \cite{Katz}).
\end{remark}

In the proposition below we show that a certain {\em stronger version}
of Theorem \ref{T:3} is false.  This answers a question of Daqing Wan,
proposed to me via email.

\begin{proposition}\label{P:counterexample}
There does not exist any Zariski dense open subset $\cU$ over $\bar\Q$
of $\A^d$ such that the following is satisfied:

For any strictly increasing sequence $\{p_i\}_{i\geq 1}$ of primes,
and for any sequence $\{f_i(x)\}_{i\geq 1}\in\cU(\bar\Q)$, where
$\cP_i$ is a prime idea of $\bar\Q(f_i)$ lying over $p_i$, one has
$$\lim_{i\rightarrow\infty}\NP(f_i\bmod\cP_i)=\HP(\A^d).$$
\end{proposition}
\begin{proof}
  Suppose there is such a Zariski dense open subset $\cU$ defined over
  $\bar\Q$ in $\A^d$.  Then $\cU$ contains the complementary set of
  zeros of $h(\vec{t})$ for some nonzero polynomial
  $h(\vec{t})\in\bar\Q[\vec{t}]$ with $\vec{t}=(t_0,\ldots,t_{d-1})$
  as the variable.  Since $h(\vec{t})\neq 0$ we also have
  $h(p\vec{t})\neq 0$ for any prime $p$.  We will construct a
  contradiction.  Choose primes $p_i\equiv -1\bmod d$ such that the
  sequence $\{p_i\}_{i\geq 1}$ is strictly increasing.  Let
  $f_i(x)=x^d+p_ic_{i,d-1}x^{d-1}+\cdots+p_ic_{i,1}x+p_ic_{i,0}$ where
  $\vec{c_i}=(c_{i,0},\ldots,c_{i,d-1})\in\bar{\Q}^{d}$ satisfies
  $h(p_i\vec{c_i})\neq 0$.  This exists because $h(p_i\vec{t})\neq 0$.
  We observe easily that $f_i\in\cU(\bar\Q)$ and $f_i(x)\equiv
  x^d\bmod \cP_i$ for every $\cP_i$ over $p_i$ in $\Q(\vec{c_i})$.
  The latter congruence implies that $\NP(f_i\bmod\cP_i)=\NP(x^d\bmod
  p_i)$.  It is well-known that for $p_i\equiv -1\bmod d$ the
  $\NP(x^d\bmod p_i)$ is a straight line of slope $1/2$ (see
  \cite[Section 6]{Zhu:1}).  Apparently this limit is not equal to the
  Hodge polygon.  This proves the proposition.
\end{proof}


\begin{thebibliography}{99}
\bibitem{Adolphson-Sperber:1}
{\sc Alan Adolphson; Steven Sperber:} On the zeta function of a
complete intersection. {\it Ann. Sci. \'Ecole Norm. Sup.} {\bf 29}
(1996), no. 3, 287--328.

\bibitem{Adolphson-Sperber}
{\sc Alan Adolphson; Steven Sperber:}
Newton polyhedra and
the degree of the $L$-function associated to an exponential sum.
{\it Invent. Math.} {\bf 88} (1987), 555--569.

\bibitem{Bombieri}
{\sc Enrico Bombieri:}
On exponential sums in finite fields,
{\it American J. Math.} {\bf 88} (1966), 71--105.

\bibitem{Dwork:1}
{\sc Bernard Dwork:}
On the zeta function of a hypersurface.
{\it Publication Math. IHES,}
(1962), 5--68.

\bibitem{Dwork:2}
{\sc Bernard Dwork:}
On the zeta function of a hypersurface. II.
{\it Ann. of Math.}
{\bf 80} (1964), 227--299.

\bibitem{Hong:1}
{\sc Shaofang Hong:}
Newton polygons of $L$ functions associated with
exponential sums of polynomials of degree four over
finite fields.
{\it Finite Fields Appl.} {\bf 7} (2001), 205--237.

\bibitem{Hong:2}
{\sc Shaofang Hong:}
Newton polygons for $L$ functions of exponential sums
of polynomials of degree six over finite fields.
{\it J. Number Theory}
{\bf 97} (2002), 368--396.

\bibitem{Katz}
{\sc Nicholas Katz:}
Slope filtration of $F$-crystals.
{\it Ast\'erisque} {\bf 63} (1979), 113--164.
Soci\'et\'e Math\'ematique de France.

\bibitem{Katz:1}
{\sc Nicholas  Katz:}
Crystalline cohomology, Dieudonn\'e modules, and Jacobi sums,
{\it Automorphic forms, Representation theory and
Arithmetic}. Tata Institute of Fundamental Research, Bombay, 1979.
165--246.

\bibitem{Katz:2}
{\sc Nicholas Katz:}
Sommes exponentielles,
{\it Ast\'erisque}
{\bf 79},
Soci\'et\'e math\'ematique de France,
1980.

\bibitem{Kedlaya:1}
{\sc Kiran Kedlaya:}
A $p$-adic local monodromy theorem.
{\tt arXiv.math.AG/0110124},
To appear in {Annals of Math.}

\bibitem{Koblitz}
{\sc Neal Koblitz:}
$p$-adic numbers, $p$-adic analysis, and Zeta-functions,
(Second edition),
{\it Graduate Texts in Mathematics} {\bf 58}.
Springer-Verlag, 1984.

\bibitem{Manin:1}
{\sc Yuri Manin:}
The theory of commutative formal groups over
fields of finite characteristic,
{\it Russian Math. Surveys}
{\bf 18} (1963),
1--83.

\bibitem{Mazur}
{\sc Barry Mazur:}
Frobenius and the Hodge filtration,
{\it Bull. of the American Math. Soc.},
{\bf 78} (1972), 653--667.

\bibitem{SZ:3}
{\sc Jasper Scholten; Hui June Zhu:}
The first slope case of Wan's conjecture.
{\it Finite fields Appl.} {\bf 8} (2002), 414--419.

\bibitem{SZ:2}
{\sc Jasper Scholten; Hui June Zhu:}
Slope estimates of Artin-Schreier curves,
{\it Compositio Math.} {\bf 137} (2003), 275--292.

\bibitem{Sperber:1}
{\sc Steven Sperber:}
Congruence properties of the hyper-Kloosterman sum.
{\it Compositio Math.}
{\bf 40} (1980), 3--33.

\bibitem{Sperber:2}
{\sc Steven Sperber:}
On the $p$-adic theory of exponential sums,
{\it American J. Math.,}
{\bf 109} (1986), 255--296.

\bibitem{van der Geer}
{\sc Gerard van der Geer; Marcel van der Vlugt:}
On the existence of supersingular
curves of given genus,
{\it J. Reine Angew. Math.,}
{\bf 458} (1995), 53--61.

\bibitem{Wan:1}
{\sc Daqing Wan:}
Newton polygons of zeta functions and $L$ functions.
{\it Annals of Mathematics} {\bf 137} (1993), 249--293.

\bibitem{Wan:2} {\sc Daqing Wan:} Variation of $p$-adic Newton
polygons of $L$ functions for exponential sums.
Preprint available at
{\tt http://www.math.uci.edu/dwan/Overview.html}.

\bibitem{Washington}
{\sc Lawrance C. Washington:}
An introduction to cyclotomic fields.
Second edition.
Graduate textbook in mathematics, vol. {\bf 83}.
Springer, 1997.

\bibitem{Zhu:1}
{\sc Hui June Zhu:}
$p$-adic variation of $L$ functions of one variable exponential
sums. I.
{\it American J. Math.} {\bf 125} (2003), 669-690.

\bibitem{Zhu:2}
{\sc Hui June Zhu:}
$L$-functions of exponential sums over one dimensional affinoids:
Newton over Hodge. {\tt http://arXiv.org/abs/math.NT/0302085.}
\end{thebibliography}
\end{document}